\documentclass[11pt]{article}
\usepackage{amssymb,amsthm,amsmath}

\def\qed{\hfill $\Box$}

\newcommand\pf{\smallbreak\noindent \texttt{Proof}. }

\begin{document}

\newtheorem{thm}{Theorem}[section]
\newtheorem{prop}[thm]{Proposition}
\newtheorem{lem}[thm]{Lemma}
\newtheorem{cor}[thm]{Corollary}
\newtheorem{ex}[thm]{Example}
\renewcommand{\thefootnote}{*}

\title{\bf Automorphism groups of some 3-dimensional Leibniz algebras}

\author{\textbf{L.A.~Kurdachenko, O.O.~Pypka}\\
Oles Honchar Dnipro National University, Dnipro, Ukraine\\
{\small e-mail: lkurdachenko@gmail.com, sasha.pypka@gmail.com}\\
\textbf{M.M.~Semko}\\
State Tax University, Irpin, Ukraine\\
{\small e-mail: dr.mykola.semko@gmail.com}}
\date{}

\maketitle

\begin{abstract}
Let $L$ be an algebra over a field $F$ with the binary operations $+$ and $[,]$. Then $L$ is called a left Leibniz algebra if it satisfies the left Leibniz identity: $[[a,b],c]=[a,[b,c]]-[b,[a,c]]$ for all elements $a,b,c\in L$. A linear transformation $f$ of $L$ is called an endomorphism of $L$, if $f([a,b])=[f(a),f(b)]$ for all elements $a,b\in L$. A bijective endomorphism of $L$ is called an automorphism of $L$. It is easy to show that the set of all automorphisms of the Leibniz algebra is a group with respect to the operation of multiplication of automorphisms. The description of the structure of the automorphism groups of Leibniz algebras is one of the natural and important problems of the general Leibniz algebra theory. The main goal of this article is to describe the structure of the automorphism group of a certain type of nilpotent three-dimensional Leibniz algebras.
\end{abstract}



\thispagestyle{empty}

\section{Introduction.}
Let $L$ be an algebra over a field $F$ with the binary operations $+$ and $[,]$. Then $L$ is called a \textit{left Leibniz algebra} if it satisfies the left Leibniz identity:
$$[[a,b],c]=[a,[b,c]]-[b,[a,c]]$$
for all elements $a,b,c\in L$.

Leibniz algebras appeared first in the paper of A.~Blokh~\cite{BA1965}, but the term ``Leibniz algebra'' appears in the book of J.-L.~Loday~\cite{LJ1992}, and the article of J.-L.~Loday~\cite{LJ1993}. In \cite{LP1993}, J.-L.~Loday and T.~Pirashvili began the real study of the properties of Leibniz algebras. The theory of Leibniz algebras has developed very intensively in many different directions. Some of the results of this theory were presented in the book~\cite{AOR2020}. Note that Lie algebras are a partial case of Leibniz algebras. Conversely, if $L$ is a Leibniz algebra, in which $[a,a]=0$ for every element $a\in L$, then it is a Lie algebra. Thus, Lie algebras can be characterized as anticommutative Leibniz algebras. At the same time, there is a very significant difference between Lie algebras and Leibniz algebras (see, for example, survey papers \cite{CPSY2019,KKPS2017,KSeSu2020,SI2021}).

Let $L$ be a Leibniz algebra. As usual, a linear transformation $f$ of $L$ is called an \textit{endomorphism} of $L$, if $f([a,b])=[f(a),f(b)]$ for all elements $a,b\in L$. Clearly, a product of two endomorphisms of $L$ is also endomorphism, so that the set of all endomorphisms of $L$ is a semigroup by a multiplication. We note that the sum of two endomorphisms of $L$ is not necessary to be an endomorphism of $L$, so that we cannot speak about the endomorphism ring of $L$.

As usual, a bijective endomorphism of $L$ is called an \textit{automorphism} of $L$. It is not hard to show that the set $Aut_{[,]}(L)$ of all automorphisms of $L$ is a group by a multiplication (see, for example, \cite{KPS2023}).

As for other algebraic structures, the search for the structure of automorphism groups of Leibniz algebras is one of the important problems of this theory. It should be noted that automorphisms groups of Leibniz algebras have hardly been studied. It is natural to start studying automorphism groups of Leibniz algebras, the structure of which has been studied quite fully. A description of the structure of automorphism groups of infinite-dimensional cyclic Leibniz algebras was obtained in \cite{KSuY2022}, and of finite-dimensional cyclic Leibniz algebras was obtained in \cite{KPS2023}. The question naturally arises about automorphism groups of Leibniz algebras of low dimension. Unlike Lie algebras, the situation with Leibniz algebras of dimension 3 is very diverse. Leibniz algebras of dimension 3 are mostly described. Their most detailed description can be found in \cite{KPS2022}. In \cite{KPV2022}, the description of automorphism groups of Leibniz algebras with dimension 3 was started. This description is quite large. The automorphism groups of only two types of nilpotent Leibniz algebras of dimension 3 are described in \cite{KPV2022}.

This article is devoted to the description of another type of nilpotent Leibniz algebras.

\section{Some preliminary results and remarks.}
Let $L$ be a Leibniz algebra over a field $F$. Then $L$ is called \textit{abelian} if $[a,b]=0$ for every elements $a,b\in L$. In particular, an abelian Leibniz algebra is a Lie algebra.

If $A,B$ are subspaces of $L$, then $[A,B]$ will denote a subspace generated by all elements $[a,b]$ where $a\in A$, $b\in B$. A subspace $A$ of $L$ is called a \textit{subalgebra} of $L$, if $[a,b]\in A$ for every $a,b\in A$. A subalgebra $A$ of $L$ is called a \textit{left} (respectively \textit{right}) \textit{ideal} of $L$, if $[b,a]\in A$ (respectively $[a,b]\in A$) for every $a\in A$, $b\in L$. A subalgebra $A$ of $L$ is called an \textit{ideal} of $L$ (more precisely, \textit{two-sided ideal}) if it is both a left ideal and a right ideal.

Every Leibniz algebra $L$ possesses the following specific ideal. Denote by $Leib(L)$ the subspace generated by the elements $[a,a]$, $a\in L$. It is not hard to prove that $Leib(L)$ is an ideal of $L$. The ideal $Leib(L)$ is called the \textit{Leibniz kernel} of $L$. We note the following important property of the Leibniz kernel: $[[a,a],x]=[a,[a,x]]-[a,[a,x]]=0$.

The \textit{left} (respectively \textit{right}) \textit{center} $\zeta^{\textrm{left}}(L)$ (respectively $\zeta^{\textrm{right}}(L)$) of a Leibniz algebra $L$ is defined by the rule:
$$\zeta^{\textrm{left}}(L)=\{x\in L|\ [x,y]=0\ \mbox{for each element }y\in L\}$$
(respectively,
$$\zeta^{\textrm{right}}(L)=\{x\in L|\ [y,x]=0\ \mbox{for each element }y\in L\}).$$
It is not hard to prove that the left center of $L$ is an ideal, but that is not true for the right center. Moreover, $Leib(L)\leqslant\zeta^{\textrm{left}}(L)$ so that $L/\zeta^{\textrm{left}}(L)$ is a Lie algebra. The right center is a subalgebra of $L$ and, in general, the left and right centers are different (see, for example, \cite{KOP2016}).

The \textit{center} $\zeta(L)$ of $L$ is defined by the rule:
$$\zeta(L)=\{x\in L|\ [x,y]=0=[y,x]\ \mbox{for each element }y\in L\}.$$
The center is an ideal of $L$.

Now we define the \textit{upper central series}
$$\langle0\rangle=\zeta_{0}(L)\leqslant\zeta_{1}(L)\leqslant\ldots\zeta_{\alpha}(L)\leqslant\zeta_{\alpha+1}(L)\leqslant\ldots\zeta_{\eta}(L)=\zeta_{\infty}(L)$$
of a Leibniz algebra $L$ by the following rule: $\zeta_{1}(L)=\zeta(L)$ is the center of $L$, and recursively, $\zeta_{\alpha+1}(L)/\zeta_{\alpha}(L)=\zeta(L/\zeta_{\alpha}(L))$ for all ordinals $\alpha$, and $\zeta_{\lambda}(L)=\bigcup_{\mu<\lambda}\zeta_{\mu}(L)$ for the limit ordinals $\lambda$. By definition, each term of this series is an ideal of $L$.

Define the \textit{lower central series} of $L$
$$L=\gamma_{1}(L)\geqslant\gamma_{2}(L)\geqslant\ldots\gamma_{\alpha}(L)\geqslant\gamma_{\alpha+1}\geqslant\ldots\gamma_{\tau}(L)=\gamma_{\infty}(L)$$
by the following rule: $\gamma_{1}(L)=L$, $\gamma_{2}(L)=[L,L]$, and recursively $\gamma_{\alpha+1}(L)=[L,\gamma_{\alpha}(L)]$ for all ordinals $\alpha$ and $\gamma_{\lambda}(L)=\bigcap_{\mu<\lambda}\gamma_{\mu}(L)$ for the limit ordinals $\lambda$.

We say that a Leibniz algebra $L$ is \textit{nilpotent}, if there exists a positive integer $k$ such that $\gamma_{k}(L)=\langle0\rangle$. More precisely, $L$ is said to be \textit{nilpotent of nilpotency class $c$} if $\gamma_{c+1}(L)=\langle0\rangle$, but $\gamma_{c}(L)\neq\langle0\rangle$. We denote the nilpotency class of $L$ by $ncl(L)$.

Let $L$ be a Leibniz algebra over a field $F$, $M$ be non-empty subset of $L$ and $H$ be a subalgebra of $L$. Put
\begin{gather*}
Ann_{H}^{\mathrm{left}}(M)=\{a\in H|\ [a,M]=\langle0\rangle\},\\
Ann_{H}^{\mathrm{right}}(M)=\{a\in H|\ [M,a]=\langle0\rangle\}.
\end{gather*}
The subset $Ann_{H}^{\mathrm{left}}(M)$ is called the \textit{left annihilator} of $M$ in subalgebra $H$. The subset $Ann_{H}^{\mathrm{right}}(M)$ is called the \textit{right annihilator} of $M$ in subalgebra $H$. The intersection
\begin{gather*}
Ann_{H}(M)=Ann_{H}^{\mathrm{left}}(M)\cap Ann_{H}^{\mathrm{right}}(M)=\\
\{a\in H|\ [a,M]=\langle0\rangle=[M,a]\}
\end{gather*}
is called the \textit{annihilator} of $M$ in subalgebra $H$. It is not hard to see that all of these subsets are subalgebras of $L$. Moreover, if $M$ is an ideal of $L$, then $Ann_{L}(M)$ is an ideal of $L$ (see, for example, \cite{KKPS2017}).

The first type of nilpotent Leibniz algebras of dimension 3 are nilpotent Leibniz algebras of nilpotency class 3. There is only one type of such algebras:
\begin{gather*}
L_{1}=Lei_{1}(3,F)=Fa_{1}\oplus Fa_{2}\oplus Fa_{3},\ \mbox{where }[a_{1},a_{1}]=a_{2},[a_{1},a_{2}]=a_{3},\\
[a_{1},a_{3}]=[a_{2},a_{1}]=[a_{2},a_{2}]=[a_{2},a_{3}]=[a_{3},a_{1}]=[a_{3},a_{2}]=[a_{3},a_{3}]=0.
\end{gather*}
This is a cyclic Leibniz algebra, $Leib(L_{1})=\zeta^{\mathrm{left}}(L_{1})=[L_{1},L_{1}]=Fa_{2}\oplus Fa_{3}$, $\zeta^{right}(L_{1})=\zeta(L_{1})=\gamma_{3}(L_{1})=Fa_{3}$.

Let now $L$ be a nilpotent Leibniz algebra, whose nilpotency class is 2. Of course, we assume that $L$ is not a Lie algebra. Then the center $\zeta(L)$ has dimension 2 or 1. Suppose that $dim_{F}(\zeta(L))=2$. Since $L$ is not a Lie algebra, there is an element $a_{1}$ such that $[a_{1},a_{1}]=a_{3}\neq0$. Since a Leibniz algebra of dimension 1 is abelian, $a_{3}\in\zeta(L)$. It follows that $[a_{1},a_{3}]=[a_{3},a_{1}]=[a_{3},a_{3}]=0$. Being an abelian algebra of dimension 2, $\zeta(L)$ has a direct decomposition $\zeta(L)=Fa_{2}\oplus Fa_{3}$ for some element $a_{2}$. Thus we come to the following type of nilpotent Leibniz algebra:
\begin{gather*}
L_{2}=Lei_{2}(3,F)=Fa_{1}\oplus Fa_{2}\oplus Fa_{3},\ \mbox{where }[a_{1},a_{1}]=a_{3},[a_{1},a_{2}]=\\
[a_{1},a_{3}]=[a_{2},a_{1}]=[a_{2},a_{2}]=[a_{2},a_{3}]=[a_{3},a_{1}]=[a_{3},a_{2}]=[a_{3},a_{3}]=0.
\end{gather*}
In other words, $L_{2}$ is a direct sum of two ideals $A=Fa_{1}\oplus Fa_{3}$ and $B=Fa_{2}$. Moreover, $A$ is a nilpotent cyclic Leibniz algebra of dimension 2, $Leib(L_{2})=[L_{2},L_{2}]=Fa_{3}$, $\zeta^{\mathrm{left}}(L_{2})=\zeta^{\mathrm{right}}(L_{2})=\zeta(L_{2})=Fa_{2}\oplus Fa_{3}$.

We note that the structure of the automorphism groups of Leibniz algebras $Lei_{1}(3,F)$ and $Lei_{2}(3,F)$ was described in \cite{KPV2022}.

Suppose now that $L$ is a nilpotent Leibniz algebra, $ncl(L)=2$ and $dim_{F}(\zeta(L))=1$. Since $L$ is not a Lie algebra, there is an element $a_{1}$ such that $[a_{1},a_{1}]=a_{3}\neq0$. Since the factor-algebra $L/\zeta(L)$ is abelian, $a_{3}\in\zeta(L)$. It follows that $[a_{1},a_{3}]=[a_{3},a_{1}]=[a_{3},a_{3}]=0$. Then $\zeta(L)=Fa_{3}$. For every element $x\in L$ we have $[a_{1},x],[x,a_{1}]\in\zeta(L)\leqslant\langle a_{1}\rangle=Fa_{1}\oplus Fa_{3}$. It follows that subalgebra $\langle a_{1}\rangle$ is an ideal of $L$. Since $dim_{F}(\langle a_{1}\rangle)=2$, $\langle a_{1}\rangle\neq L$. Choose an element $b$ such that $b\not\in\langle a_{1}\rangle$. We have $[b,a_{1}]=\gamma a_{3}$ for some $\gamma\in F$. If $\gamma\neq0$, then put $b_{1}=\gamma^{-1}b-a_{1}$. Then $[b_{1},a_{1}]=0$. The choice of $b_{1}$ shows that $b_{1}\not\in\langle a_{1}\rangle$. If follows that subalgebra $Ann^{\mathrm{left}}_{L}(a_{1})$ has dimension 2. Suppose first that $Ann^{\mathrm{left}}_{L}(a_{1})$ is an abelian subalgebra. Then it has a direct decomposition $Ann^{\mathrm{left}}_{L}(a_{1})=Fa_{2}\oplus Fb_{2}$ for some element $b_{2}$, where $[b_{2},b_{2}]=0$. Since $dim_{F}(\zeta(L))=1$, $b_{2}\not\in\zeta(L)$. Then $[a_{1},b_{2}]=\lambda a_{3}$ where $0\neq\lambda\in F$. Put now $a_{2}=\lambda^{-1}b_{2}$. Thus, we come to the following type of nilpotent Leibniz algebra:
\begin{gather*}
L_{3}=Lei_{3}(3,F)=Fa_{1}\oplus Fa_{2}\oplus Fa_{3},\ \mbox{where }[a_{1},a_{1}]=[a_{1},a_{2}]=a_{3},\\
[a_{1},a_{3}]=[a_{2},a_{1}]=[a_{2},a_{2}]=[a_{2},a_{3}]=[a_{3},a_{1}]=[a_{3},a_{2}]=[a_{3},a_{3}]=0.
\end{gather*}
In other words, $L_{3}$ is a direct sum of ideal $A=Fa_{1}\oplus Fa_{3}$ and subalgebra $B=Fa_{2}$. Moreover, $A$ is a nilpotent cyclic Leibniz algebra of dimension 2, $Leib(L_{3})=[L_{3},L_{3}]=\zeta^{\mathrm{right}}(L_{3})=\zeta(L_{3})=Fa_{3}$, $\zeta^{\mathrm{left}}(L_{3})=Fa_{2}\oplus Fa_{3}$.

This article is devoted to the description of this type of nilpotent Leibniz algebras.

Here are some general useful properties of automorphism groups of Leibniz algebras. Their proofs are given in the article \cite{KPV2022}.

\begin{lem}\label{L1}
Let $L$ be a Leibniz algebra over a field $F$ and $f$ be an automorphism of $L$. Then $f(\zeta^{\mathrm{left}}(L))=\zeta^{\mathrm{left}}(L)$, $f(\zeta^{\mathrm{right}}(L))=\zeta^{\mathrm{right}}(L)$, $f(\zeta(L))=\zeta(L)$, $f([L,L])=[L,L]$.
\end{lem}

\begin{lem}\label{L2}
Let $L$ be a Leibniz algebra over a field $F$ and $f$ be an automorphism of $L$. Then $f(\zeta_{\alpha}(L))=\zeta_{\alpha}(L)$, $f(\gamma_{\alpha}(L))=\gamma_{\alpha}(L)$ for all ordinals $\alpha$. In particular, $f(\zeta_{\infty}(L))=\zeta_{\infty}(L)$ and $f(\gamma_{\infty}(L))=\gamma_{\infty}(L)$.
\end{lem}

\begin{lem}\label{L3}
Let $L$ be a Leibniz algebra over a field $F$ and $f$ be an endomorphism of $L$. Then $f(\gamma_{\alpha}(L))\leqslant\gamma_{\alpha}(L)$ for all ordinals $\alpha$. In particular, $f(\gamma_{\infty}(L))\leqslant\gamma_{\infty}(L)$.
\end{lem}

Let $L$ be a Leibniz algebra over a field $F$, $A$ be a subalgebra of $L$, $G=Aut_{[,]}(L)$. Put
$$C_{G}(A)=\{\alpha\in G|\ \alpha(x)=x\ \mbox{for every element }x\in A\}.$$
If $A$ is an ideal of $L$, then put
$$C_{G}(L/A)=\{\alpha\in G|\ \alpha(x+A)=x+A\ \mbox{for every element }x\in L\}.$$

\begin{lem}\label{L4}
Let $L$ be a Leibniz algebra over a field $F$ and $G=Aut_{[,]}(L)$. If $A$ is a $G$-invariant subalgebra, then $C_{G}(A)$ and $C_{G}(L/A)$ are normal subgroups of $G$.
\end{lem}

\section{Main result.}
\begin{thm}\label{T1}
Let $G$ be an automorphism group of Leibniz algebra $L_{3}$. Then $G$ is isomorphic to a subgroup of $GL_{3}(F)$, consisting of the matrices, having the following form:
\begin{equation*}
\left(\begin{array}{ccc}
\alpha_{1} & 0 & 0\\
\alpha_{2} & \alpha_{1}+\alpha_{2} & 0\\
\alpha_{3} & \beta_{3} & \alpha_{1}^{2}+\alpha_{1}\alpha_{2}
\end{array}\right)
\end{equation*}
where $\alpha_{1}\neq0$, $\alpha_{1}+\alpha_{2}\neq0$. This subgroup is a semidirect product of normal subgroup $S_{3}(L,F)$, which is isomorphic to a subgroup of $GL_{3}(F)$, consisting of the matrices of the form
\begin{equation*}
\left(\begin{array}{ccc}
1 & 0 & 0\\
\alpha_{2} & 1+\alpha_{2} & 0\\
\alpha_{3} & \beta_{3} & 1+\alpha_{2}
\end{array}\right)
\end{equation*}
and a subgroup $D_{3}(L,F)$, consisting of the matrices of the form
\begin{equation*}
\left(\begin{array}{ccc}
\sigma & 0 & 0\\
0 & \sigma & 0\\
0 & 0 & \sigma^{2}
\end{array}\right).
\end{equation*}
In particular, $D_{3}(L,F)$ is isomorphic to multiplicative group of a field $F$. Furthermore, $S_{3}(L,F)$ is a semidirect product of subgroup $T_{3}(L,F)$, which is normal in $G$ and isomorphic to a subgroup of $GL_{3}(F)$, consisting of the matrices of the form
\begin{equation*}
\left(\begin{array}{ccc}
1 & 0 & 0\\
0 & 1 & 0\\
\alpha_{3} & \beta_{3} & 1
\end{array}\right),
\end{equation*}
and a subgroup $J_{3}(L,F)$, which is isomorphic to a subgroup of $GL_{3}(F)$, consisting of the matrices of the form
\begin{equation*}
\left(\begin{array}{ccc}
1 & 0 & 0\\
\lambda & 1+\lambda & 0\\
0 & 0 & 1+\lambda
\end{array}\right).
\end{equation*}
A subgroup $T_{3}(L,F)$ is isomorphic to direct product of two copy of additive group of a field $F$, and a subgroup $J_{3}(L,F)$ is isomorphic to multiplicative group of a field $F$.
\end{thm}
\pf Let $L=Lei_{3}(3,F)$, $f\in Aut_{[,]}(L)$. By Lemma~\ref{L1}, $f(Fa_{3})=Fa_{3}$, $f(Fa_{2}\oplus Fa_{3})=Fa_{2}\oplus Fa_{3}$, so that
\begin{gather*}
f(a_{1})=\alpha_{1}a_{1}+\alpha_{2}a_{2}+\alpha_{3}a_{3},\\
f(a_{2})=\beta_{2}a_{2}+\beta_{3}a_{3}.
\end{gather*}
Then
\begin{gather*}
f(a_{3})=f([a_{1},a_{1}])=[f(a_{1}),f(a_{1})]=\\
[\alpha_{1}a_{1}+\alpha_{2}a_{2}+\alpha_{3}a_{3},\alpha_{1}a_{1}+\alpha_{2}a_{2}+\alpha_{3}a_{3}]=\\
\alpha_{1}^{2}[a_{1},a_{1}]+\alpha_{1}\alpha_{2}[a_{1},a_{2}]=\alpha_{1}^{2}a_{3}+\alpha_{1}\alpha_{2}a_{3}=(\alpha_{1}^{2}+\alpha_{1}\alpha_{2})a_{3};\\
f(a_{3})=f([a_{1},a_{2}])=[f(a_{1}),f(a_{2})]=[\alpha_{1}a_{1}+\alpha_{2}a_{2}+\alpha_{3}a_{3},\beta_{2}a_{2}+\beta_{3}a_{3}]=\\
\alpha_{1}\beta_{2}[a_{1},a_{2}]=\alpha_{1}\beta_{2}a_{3}.
\end{gather*}
Thus, we obtain an equality $\alpha_{1}(\alpha_{1}+\alpha_{2})=\alpha_{1}\beta_{2}$. Being an automorphism, $f$ is a non-degenerate linear transformation, so that $\alpha_{1}\neq0$. It follows that $\alpha_{1}+\alpha_{2}=\beta_{2}$. Thus, an automorphism $f$ has in basis $\{a_{1},a_{2},a_{3}\}$ the following matrix
\begin{equation*}
\left(\begin{array}{ccc}
\alpha_{1} & 0 & 0\\
\alpha_{2} & \alpha_{1}+\alpha_{2} & 0\\
\alpha_{3} & \beta_{3} & \alpha_{1}^{2}+\alpha_{1}\alpha_{2}
\end{array}\right)
\end{equation*}
Denote by $\Xi$ the canonical monomorphism of $End_{[,]}(L)$ in $M_{3}(F)$.

Put
$$S=\{f|\ f\in End(L),f(a_{1})\in a_{1}+\zeta^{\mathrm{left}}(L)\}=C_{End(L)}(L/\zeta^{\mathrm{left}}(L)).$$
If $f\in S\cap Aut_{[,]}(L)$, then $f(a_{1})=a_{1}+\alpha_{2}a_{2}+\alpha_{3}a_{3}$, $f(a_{2})=(1+\alpha_{2})a_{2}+\beta_{3}a_{3}$, $f(a_{3})=(1+\alpha_{2})a_{3}$. If $x$ is an arbitrary element of $L$, $x=\xi_{1}a_{1}+\xi_{2}a_{2}+\xi_{3}a_{3}$, then
\begin{gather*}
f(x)=\xi_{1}f(a_{1})+\xi_{2}f(a_{2})+\xi_{3}f(a_{3})=\\
\xi_{1}a_{1}+\xi_{1}\alpha_{2}a_{2}+\xi_{1}\alpha_{3}a_{3}+\xi_{2}((1+\alpha_{2})a_{2}+\beta_{3}a_{3})+\xi_{3}(1+\alpha_{2})a_{3}=\\
\xi_{1}a_{1}+(\xi_{1}\alpha_{2}+\xi_{2}+\xi_{2}\alpha_{2})a_{2}+(\xi_{1}\alpha_{3}+\xi_{2}\beta_{3}+\xi_{3}(1+\alpha_{2}))a_{3}.
\end{gather*}

Conversely, let $\lambda,\mu,\nu$ be the elements of $F$, $v_{\lambda,\mu,\nu}$ be a linear transformation of $L$, defined by the rule: if $x=\xi_{1}a_{1}+\xi_{2}a_{2}+\xi_{3}a_{3}$, then
$$v_{\lambda,\mu,\nu}(x)=\xi_{1}a_{1}+(\xi_{1}\lambda+\xi_{2}+\xi_{2}\lambda)a_{2}+(\xi_{1}\mu+\xi_{2}\nu+\xi_{3}(1+\lambda))a_{3}.$$
Let $x,y$ be the arbitrary elements of $L$, $x=\xi_{1}a_{1}+\xi_{2}a_{2}+\xi_{3}a_{3}$, $y=\eta_{1}a_{1}+\eta_{2}a_{2}+\eta_{3}a_{3}$, where $\xi_{1},\xi_{2},\xi_{3},\eta_{1},\eta_{2},\eta_{3}\in F$. Then
\begin{gather*}
[x,y]=[\xi_{1}a_{1}+\xi_{2}a_{2}+\xi_{3}a_{3},\eta_{1}a_{1}+\eta_{2}a_{2}+\eta_{3}a_{3}]=\\
\xi_{1}\eta_{1}[a_{1},a_{1}]+\xi_{1}\eta_{2}[a_{1},a_{2}]=\xi_{1}(\eta_{1}+\eta_{2})a_{3};\\
v_{\lambda,\mu,\nu}([x,y])=v_{\lambda,\mu,\nu}(\xi_{1}(\eta_{1}+\eta_{2})a_{3})=\xi_{1}(\eta_{1}+\eta_{2})v_{\lambda,\mu,\nu}(a_{3})=\\
\xi_{1}(\eta_{1}+\eta_{2})(1+\lambda)a_{3};\\
[v_{\lambda,\mu,\nu}(x),v_{\lambda,\mu,\nu}(y)]=\\
[\xi_{1}a_{1}+(\xi_{1}\lambda+\xi_{2}+\xi_{2}\lambda)a_{2}+(\xi_{1}\mu+\xi_{2}\nu+\xi_{3}(1+\lambda))a_{3},\\
\eta_{1}a_{1}+(\eta_{1}\lambda+\eta_{2}+\eta_{2}\lambda)a_{2}+(\eta_{1}\mu+\eta_{2}\nu+\eta_{3}(1+\lambda))a_{3}]=\\
\xi_{1}\eta_{1}[a_{1},a_{1}]+\xi_{1}(\eta_{1}\lambda+\eta_{2}+\eta_{2}\lambda)[a_{1},a_{2}]=(\xi_{1}\eta_{1}+\xi_{1}(\eta_{1}\lambda+\eta_{2}+\eta_{2}\lambda))a_{3}=\\
\xi_{1}(\eta_{1}+\eta_{1}\lambda+\eta_{2}+\eta_{2}\lambda)a_{3}=\xi_{1}(\eta_{1}+\eta_{2})(1+\lambda)a_{3},
\end{gather*}
so that $v_{\lambda,\mu,\nu}([x,y])=[v_{\lambda,\mu,\nu}(x),v_{\lambda,\mu,\nu}(y)]$. It shows that $S\leqslant Aut_{[,]}(L)$. Moreover, by Lemma~\ref{L4}, $S$ is a normal subgroup of $Aut_{[,]}(L)$. Furthermore, put $S_{3}(L,F)=\Xi(S)$. Then $S_{3}(L,F)$ is a subgroup of a group $T_{3}(F)$, which consist of the matrices, having the following form:
\begin{equation*}
\left(\begin{array}{ccc}
1 & 0 & 0\\
\alpha_{2} & 1+\alpha_{2} & 0\\
\alpha_{3} & \beta_{3} & 1+\alpha_{2}
\end{array}\right).
\end{equation*}

Let
\begin{gather*}
T=\{f|\ f\in End(L),f(a_{1})\in a_{1}+[L,L],f(a_{2})\in a_{2}+[L,L]\}=\\
C_{End(L)}(L/[L,L]).
\end{gather*}
If $f\in T\cap Aut_{[,]}(L)$, then $f(a_{1})=a_{1}+\alpha_{3}a_{3}$, $f(a_{2})=a_{2}+\beta_{3}a_{3}$, $f(a_{3})=a_{3}$. If $x$ is an arbitrary element of $L$, $x=\xi_{1}a_{1}+\xi_{2}a_{2}+\xi_{3}a_{3}$, then
\begin{gather*}
f(x)=\xi_{1}f(a_{1})+\xi_{2}f(a_{2})+\xi_{3}f(a_{3})=\\
\xi_{1}a_{1}+\xi_{1}\alpha_{3}a_{3}+\xi_{2}a_{2}+\xi_{2}\beta_{3}a_{3}+\xi_{3}a_{3}=\\
\xi_{1}a_{1}+\xi_{2}a_{2}+(\xi_{1}\alpha_{3}+\xi_{2}\beta_{3}+\xi_{3})a_{3}.
\end{gather*}
Conversely, let $\lambda,\mu$ be the elements of $F$, $z_{\lambda,\mu}$ be a linear transformation of $L$, defined by the rule: if $x=\xi_{1}a_{1}+\xi_{2}a_{2}+\xi_{3}a_{3}$, then
$$z_{\lambda,\mu}(x)=\xi_{1}a_{1}+\xi_{2}a_{2}+(\xi_{1}\lambda+\xi_{2}\mu+\xi_{3})a_{3}.$$
Let $x,y$ be the arbitrary elements of $L$, $x=\xi_{1}a_{1}+\xi_{2}a_{2}+\xi_{3}a_{3}$, $y=\eta_{1}a_{1}+\eta_{2}a_{2}+\eta_{3}a_{3}$, where $\xi_{1},\xi_{2},\xi_{3},\eta_{1},\eta_{2},\eta_{3}\in F$. Then
\begin{gather*}
[x,y]=[\xi_{1}a_{1}+\xi_{2}a_{2}+\xi_{3}a_{3},\eta_{1}a_{1}+\eta_{2}a_{2}+\eta_{3}a_{3}]=\xi_{1}(\eta_{1}+\eta_{2})a_{3};\\
z_{\lambda,\mu}([x,y])=z_{\lambda,\mu}(\xi_{1}(\eta_{1}+\eta_{2})a_{3})=\xi_{1}(\eta_{1}+\eta_{2})z_{\lambda,\mu}(a_{3})=\xi_{1}(\eta_{1}+\eta_{2})a_{3};\\
[z_{\lambda,\mu}(x),z_{\lambda,\mu}(y)]=\\
[\xi_{1}a_{1}+\xi_{2}a_{2}+(\xi_{1}\lambda+\xi_{2}\mu+\xi_{3})a_{3},\eta_{1}a_{1}+\eta_{2}a_{2}+(\eta_{1}\lambda+\eta_{2}\mu+\eta_{3})a_{3}]=\\
\xi_{1}\eta_{1}[a_{1},a_{1}]+\xi_{1}\eta_{2}[a_{1},a_{2}]=\xi_{1}(\eta_{1}+\eta_{2})a_{3},
\end{gather*}
so that $z_{\lambda,\mu}([x,y])=[z_{\lambda,\mu}(x),z_{\lambda,\mu}(y)]$. It shows that $T$ is a subgroup of $Aut_{[,]}(L)$. Furthermore, put $T_{3}(L,F)=\Xi(T)$. Then $T_{3}(L,F)$ is a subgroup of a group $UT_{3}(F)$ of all unitriangular matrices over a field $F$, which consist of the matrices, having the following form:
\begin{equation*}
\left(\begin{array}{ccc}
1 & 0 & 0\\
0 & 1 & 0\\
\alpha_{3} & \beta_{3} & 1
\end{array}\right).
\end{equation*}
It is not hard to see, that $T_{3}(L,F)$ is abelian and it is isomorphic to direct product of two copy of additive group of a field $F$. Clearly $\Xi(Aut_{[,]}(L))\cap UT_{3}(F)=T_{3}(L,F)$, so it follows that a subgroup $T$ is normal in $Aut_{[,]}(L)$.

Let
$$J=\{f|\ f\in S,f(a_{1})=a_{1}+\lambda a_{2},f(a_{2})=(1+\lambda)a_{2},f(a_{3})=(1+\lambda)a_{3},\lambda\in F\}.$$
Put $J_{3}(L,F)=\Xi(J)$. Then $J_{3}(L,F)$ is a subset of $T_{3}(F)$, which consist of the matrices, having the following form:
\begin{equation*}
\left(\begin{array}{ccc}
1 & 0 & 0\\
\lambda & 1+\lambda & 0\\
0 & 0 & 1+\lambda
\end{array}\right).
\end{equation*}
An equality
\begin{gather*}
\left(\begin{array}{ccc}
1 & 0 & 0\\
\lambda & 1+\lambda & 0\\
0 & 0 & 1+\lambda
\end{array}\right)
\left(\begin{array}{ccc}
1 & 0 & 0\\
\mu & 1+\mu & 0\\
0 & 0 & 1+\mu
\end{array}\right)=\\
\left(\begin{array}{ccc}
1 & 0 & 0\\
(1+\lambda)(1+\mu)-1 & (1+\lambda)(1+\mu) & 0\\
0 & 0 & (1+\lambda)(1+\mu)
\end{array}\right)
\end{gather*}
shows that $J_{3}(L,F)$ is a subgroup of $S_{3}(L,F)$. Moreover, it is not hard to see that $J_{3}(L,F)$ is isomorphic to multiplicative group of a field $F$. Furthermore, it is not hard to see that the matrix equation
\begin{gather*}
\left(\begin{array}{ccc}
1 & 0 & 0\\
\alpha_{2} & 1+\alpha_{2} & 0\\
\alpha_{3} & \beta_{3} & 1+\alpha_{2}
\end{array}\right)=
\left(\begin{array}{ccc}
1 & 0 & 0\\
0 & 1 & 0\\
x & y & 1
\end{array}\right)
\left(\begin{array}{ccc}
1 & 0 & 0\\
z & 1+z & 0\\
0 & 0 & 1+z
\end{array}\right)=\\
\left(\begin{array}{ccc}
1 & 0 & 0\\
z & 1+z & 0\\
x+yz & y+yz & 1+z
\end{array}\right)
\end{gather*}
has the solutions. It proves that $S_{3}(L,F)=T_{3}(L,F)J_{3}(L,F)$ and therefore $S=TJ$. Clearly the intersection $T\cap J$ is trivial.

Let
$$D=\{f|\ f\in Aut_{[,]}(L),f(a_{1})=\sigma a_{1},f(a_{2})=\sigma a_{2},\sigma\in F\}.$$
By proved above, $f(a_{3})=\sigma^{2}a_{3}$.

Put $D_{3}(L,F)=\Xi(D)$. Then $D_{3}(L,F)$ is a subset of $T_{3}(F)$, which consist of the matrices, having the following form:
\begin{equation*}
\left(\begin{array}{ccc}
\sigma & 0 & 0\\
0 & \sigma & 0\\
0 & 0 & \sigma^{2}
\end{array}\right).
\end{equation*}
An equality
\begin{equation*}
\left(\begin{array}{ccc}
\sigma & 0 & 0\\
0 & \sigma & 0\\
0 & 0 & \sigma^{2}
\end{array}\right)
\left(\begin{array}{ccc}
\nu & 0 & 0\\
0 & \nu & 0\\
0 & 0 & \nu^{2}
\end{array}\right)=
\left(\begin{array}{ccc}
\sigma\nu & 0 & 0\\
0 & \sigma\nu & 0\\
0 & 0 & \sigma^{2}\nu^{2}
\end{array}\right)
\end{equation*}
shows that $D_{3}(L,F)$ is a subgroup of $\Xi(Aut_{[,]}(L))$. Moreover, it is not hard to see that $D_{3}(L,F)$ is isomorphic to multiplicative group of a field $F$. Furthermore, an equality
\begin{gather*}
\left(\begin{array}{ccc}
\alpha_{1} & 0 & 0\\
\alpha_{2} & \alpha_{1}+\alpha_{2} & 0\\
\alpha_{3} & \beta_{3} & \alpha_{1}^{2}+\alpha_{1}\alpha_{2}
\end{array}\right)=\\
\left(\begin{array}{ccc}
1 & 0 & 0\\
\alpha_{2}\alpha_{1}^{-1} & 1+\alpha_{2}\alpha_{1}^{-1} & 0\\
\alpha_{3}\alpha_{1}^{-1} & \beta_{3}\alpha_{1}^{-1} & 1+\alpha_{2}\alpha_{1}^{-1}
\end{array}\right)
\left(\begin{array}{ccc}
\alpha_{1} & 0 & 0\\
0 & \alpha_{1} & 0\\
0 & 0 & \alpha_{1}^{2}
\end{array}\right)
\end{gather*}
proves that $\Xi(Aut_{[,]}(L))=S_{3}(L,F)D_{3}(L,F)$ and therefore $Aut_{[,]}(L)=SD$. Clearly the intersection $S\cap D$ is trivial. Thus, we obtain that $Aut_{[,]}(L)=S\leftthreetimes D$, $D\cong F^{\times}$, $S=T\leftthreetimes J$, moreover $T$ is normal in $Aut_{[,]}(L)$, $T\cong F_{+}\times F_{+}$, $J\cong F^{\times}$. \qed

\end{document}